\documentclass{article}

\usepackage[a4paper]{geometry}

\usepackage{amsmath}

\usepackage{amssymb}
\usepackage{theorem}
\usepackage{enumerate}

\newcommand{\ov}{\overline}
\newcommand{\C}{\mathbb C}
\newcommand{\F}{\mathbb F}

\newcommand{\R}{\mathbb R}
\newcommand{\Z}{\mathbb Z}

\newcommand{\Aut}{\operatorname{Aut}}

\newcommand{\tr}{\operatorname{tr}}

\newcommand{\im}{\operatorname{Im}}

\newcommand{\eop}{\hspace*{\fill} $\Box$}

\newcommand{\SL}{\operatorname{SL}}
\newcommand{\Sp}{\operatorname{Sp}}
\newcommand{\GL}{\operatorname{GL}}
\newcommand{\Or}{\operatorname{O}}
\newcommand{\A}{\mathcal{A}}

\newcommand{\M}{\operatorname{M}}
\newcommand{\Stab}{\operatorname{Stab}}
\renewcommand{\S}{\operatorname{S}}

\theoremstyle{break}

\newtheorem{thm}{Theorem}[section]
\newtheorem{prp}[thm]{Proposition}
\newtheorem{cor}[thm]{Corollary}
\newtheorem{lem}[thm]{Lemma}
\begin{document}

\begin{center}
{\Large \bf Harmonic theta series \\[2mm]
and the Kodaira dimension of $\mathcal{A}_6$}\\[10mm]
%Moritz Dittmann\textsuperscript{1}, Riccardo Salvati Manni\textsuperscript{2}, Nils R.\ Scheithauer\textsuperscript{1}\\[4mm]
%\textsuperscript{1}{Fachbereich Mathematik, Technische Universit\"{a}t Darmstadt, Darmstadt, Germany}\\
%\textsuperscript{2}{Universit{\'a} di Roma "La Sapienza", Roma, Italy}
%\textsuperscript{2}{Dipartimento di Matematica "Guido Castelnuovo", Universit{\'a} di Roma "La  Sapienza", Italy} 
Moritz Dittmann\footnote[1]{Fachbereich Mathematik, Technische Universit\"{a}t Darmstadt, Darmstadt, Germany},
Riccardo Salvati Manni\footnote[2]{Dipartimento di Matematica "Guido Castelnuovo", Universit{\`a} di Roma "La  Sapienza", Italy},
Nils R.\ Scheithauer\footnotemark[1]
\end{center}

\vspace*{10mm}
\noindent
We construct a basis of the space $\S_{14}(\Sp_{12}(\Z))$ of Siegel cusp forms of degree $6$ and weight $14$ consisting of harmonic theta series. One of these functions has vanishing order $2$ at the boundary which implies that the Kodaira dimension of $\A_6$ is non-negative.

\vspace*{12mm}
\section{Introduction}

Siegel modular forms are natural generalizations of elliptic modular forms on $\SL_2(\Z)$. They play an important role in number theory and algebraic geometry. One of the main problems in the theory is that it is very difficult to determine the Fourier coefficients of these forms even in simple examples. In this paper we construct a basis of the space $\S_{14}(\Sp_{12}(\Z))$ consisting of harmonic theta series for Niemeier lattices and calculate some of their Fourier coefficients (see Theorems \ref{coefficients} and \ref{basis}). As a consistency check we also compute the Fourier coefficients of the Ikeda lift of the unique normalized cusp form in $\S_{22}(\SL_{2}(\Z))$.   

The quotient $\A_n = \Sp_{2n}(\Z) \backslash H_n$ is the (coarse) moduli space of principally polarized complex abelian varieties of dimension $n$. 
It was believed for a long time that the spaces $\A_n$ are unirational. This was disproved by Freitag for $n = 1 \! \mod 8$, $n \geq 17$ \cite{F1} and for $n = 0 \! \mod 24$ \cite{F2}. In the latter work he also showed that the Kodaira dimension of $\A_n$ tends to infinity for $n = 0 \! \mod 24$. A few years later Tai proved that $\A_n$ is of general type for all $n \geq 9$ \cite{Tai}. This was refined by Freitag to $n \geq 8$ \cite{F} and by Mumford to $n \geq 7$ \cite{Mum}. The ideas of Freitag and Tai strongly influenced the classification of the moduli spaces ${\mathcal M}_n$ of complex curves of genus $n$ \cite{HM}. On the other hand, $\A_n$ is unirational for $n \leq 5$. This is a classical result for $n \leq 3$ and was proved for $n = 4,5$ by Clemens \cite{C}, Donagi \cite{D}, Mori and Mukai \cite{MoMu} and Verra \cite{V}. Hence $\A_n$ is of general type for $n \geq 7$ and has Kodaira dimension $-\infty$ for $n \leq 5$. The case $n=6$ has been open for more than 30 years. One of the basis elements that we construct is a harmonic theta series related to the Leech lattice. Since the Leech lattice contains no vectors of norm 2, this function defines a global section of $\omega_{\ov{\A_6}}^2$. This implies that the Kodaira dimension of $\A_6$ is non-negative (see Theorem \ref{kodairadim}). We remark that this result does not require the knowledge of a full basis of $\S_{14}(\Sp_{12}(\Z))$. 

The paper is organized as follows. In section 2 we recall some standard results on Siegel modular forms. In the next section we describe a method to compute the Fourier coefficients of harmonic theta series. In section 4 we construct a basis of $\S_{14}(\Sp_{12}(\Z))$ consisting of harmonic theta series $\theta_{N,h,2}$. Up to a scalar there is a unique cusp form with vanishing order 2 at the boundary. We use this function in section 5 to show that the Kodaira dimension of $\mathcal{A}_6$ is non-negative. In the appendix we list the lattices $N$ and the matrices $h$ of the $\theta_{N,h,2}$.\\[-2mm]

The authors thank E.\ Freitag, M.\ M\"oller and A.\ Verra for stimulating discussions and the referees for suggesting several improvements. They also acknowledge support from the LOEWE research unit \emph{Uniformized Structures in Arithmetic and Geometry}.

\section{Siegel modular forms} 

In this section we recall some results on Siegel modular forms from \cite{F}, \cite{dimformula}, \cite{Basisproblem}, \cite{Ikedalift} and \cite{Likeda}.\\[-2mm]

Let $E \in \M_n(\C)$ be the identity matrix and $\Omega = \left( \begin{smallmatrix} 0 & E  \\ -E & 0 \end{smallmatrix} \right)$. The symplectic group
\[  \Sp_{2n}(\Z) = \{ M \in \GL_{2n}(\Z) \, | \, M^T \Omega M = \Omega \}  \]
acts on the Siegel upper halfspace 
\[  H_n = \{ Z \in \M_n(\C) \, | \, Z = Z^T, \im(Z) >0 \}  \]
by
\[  MZ = (AZ+B)(CZ+D)^{-1}  \]
where $M = \left( \begin{smallmatrix} A & B \\ C & D \end{smallmatrix} \right) \in \Sp_{2n}(\Z)$ and $Z \in H_n$. 

Let $k$ be an integer. 
%$k\in \Z$. 
A holomorphic function $f \! : H_n \to \C$ is called a Siegel modular form of degree $n$ and weight $k$ if
\[ f(MZ) = \det(CZ+D)^kf(Z) \]
for all $M = \left( \begin{smallmatrix} A & B \\ C & D \end{smallmatrix} \right) \in \Sp_{2n}(\Z)$ and $f$ is bounded in the domain $\{ Z \in H_n \, | \, \im(Z) - Y \geq 0 \}$ for all positive definite $Y \in \M_n(\R)$. For $n>1$ the second condition can be dropped by the Koecher principle.

Let $T \in \M_n(\Z)$ be symmetric. $T$ is called even if the diagonal entries of $T$ are even. If $T$ is positive semidefinite, we define 
\[  m(T) = \min_{x \in \Z^n \backslash \{0 \}} \frac{1}{2} x^TTx  \, . \]
   
Let $f$ be a Siegel modular form of degree $n$ and weight $k$. Then $f$ has a Fourier expansion of the form
\[  f(Z) = \sum_{\substack{T \in \M_n(\Z) \\ T \geq 0 \text{ and even} }} 
                                               a(T) e^{\pi i \tr(TZ)}  \]
with Fourier coefficients $a(T) \in \C$ satisfying 
\[  a(U^TTU) = \det(U)^k a(T)  \]
for all $U \in \GL_n(\Z)$. The function $f$ is called a cusp form if the vanishing order 
\[  m(f) = 
\min_{\substack{T \in \M_n(\Z) \\ T \geq 0 \text{ and even} \\ a(T) \neq 0 }}  m(T) \]
of $f$ at $\infty$ is positive, i.e. the Fourier expansion of $f$ is supported only on positive definite matrices. 

We denote the space of Siegel modular forms of degree $n$ and weight $k$ by $\M_k(\Sp_{2n}(\Z))$ and the subspace of cusp forms by $\S_k(\Sp_{2n}(\Z))$. Ta\"ibi \cite{dimformula} gives an algorithm based on Arthur's trace formula to determine the dimensions of these spaces.

An important class of Siegel modular forms is given by harmonic theta series.

Let $V$ be a euclidean vector space of dimension $m$ with scalar product $(\, , \, )$ and $L \subset V$ an even unimodular lattice of rank $m$. Later we will write the elements of $V$ as column vectors. For two $n$-tuples $x=(x_1, \ldots, x_n)$, $y=(y_1, \ldots, y_n)$ in $V^n_{\C}$ with $V_{\C} = V \otimes_{\R} \C$ we define the matrix $Q(x,y)$ of inner products, i.e.
\[ Q(x,y) =  ( (x_i,y_j) ) \, . \]
If $x=y$, we use the abbreviation $Q(x) = Q(x,x)$. Now we choose $h=(h_1, \ldots, h_n) \in V^n_{\C}$ such that 
\[ Q(h) = Q(h,h) = 0 \] 
and
\[ Q(h,\ov{h}) > 0 \, . \] 
For a positive integer $k$ define 
\[ \theta_{L,h,k}(Z) = \sum_{x \in L^n} \det(Q(x,h))^k e^{\pi i \tr(Q(x)Z) }  \, . \]
Then $\theta_{L,h,k}$ is a Siegel cusp form of degree $n$ and weight $m/2 + k$. The Fourier expansion of $\theta_{L,h,k}$ is given by 
\[ \theta_{L,h,k}(Z) = \sum_{\substack{T \in \M_n(\Z) \\ T \geq 0 \text{ and even} }} 
                                               a(T) e^{\pi i \tr(TZ)}  \]
with 
\[ a(T) = \sum_{\substack{x \in L^n \\ Q(x)=T}} \det(Q(x,h))^k \, . \]
Note that the number of summands in $a(T)$ can be very large. B\"ocherer \cite{Basisproblem} showed 

\begin{prp} \label{both}
Let $n,m$ and $k$ be positive integers with $8|m$ and $m/2>2n$. Then every Siegel cusp form of degree $n$ and weight $m/2+k$ is a linear combination of harmonic Siegel theta series of the form $\theta_{L,h,k}$ where $L$ is an even unimodular lattice of rank $m$.    
\end{prp}

\noindent
We will be interested in the case $m=24$ and $n=6$, which is not covered by this proposition.

Another construction of Siegel cusp forms is the Ikeda lift \cite{Ikedalift}. Let $k$ and $n$ be two positive integers such that $k = n \! \mod 2$ and $f \in \S_{2k}(\SL_2(\Z))$ a normalized Hecke eigenform. Then the Ikeda lift $F$ of $f$ is a Siegel cusp form of degree $2n$ and weight $k+n$ such that the standard zeta function $L(F,s)$ of $F$ is related to the $L$-function of $f$ by
\[  L(F,s) = \zeta(s) \prod_{j=1}^{2n} L(f,s+k+n-j)  \, . \]
Kohnen \cite{Likeda} gave a linear version of the Ikeda lift which can be described as follows. 

\begin{thm} \label{Ikedalift}
Let $f \in \S_{2k}(\SL_2(\Z))$ be a normalized Hecke eigenform and 
\[  g(\tau) = \sum_{\substack{m \geq 1 \\ (-1)^km = 0,1 \bmod 4}} c(m) q^m  \in 
          \S^+_{k+1/2}(\tilde{\Gamma}_0(4))                          \]
a Hecke eigenform corresponding to $f$ under the Shimura correspondence.
For a positive definite even matrix $T$ of size $2n$ write $(-1)^n \det(T) = D_T d_T^2$ where $D_T$ is a fundamental discriminant and $d_T$ a positive integer. 
Then the Ikeda lift $F \in \S_{k+n}(\Sp_{4n}(\Z))$ of $f$ is given by
\[  F(Z) = \sum_{\substack{T \in \M_{2n}(\Z) \\ T > 0 \text{ and even} }} a(T) e^{\pi i \tr(TZ)}  \]
with
\[  a(T)=\sum_{0<a|d_T} a^{k-1} \phi(a,T) 
    c \! \left( |D_T| \bigg(\frac{d_T}{a}\bigg)^2 \right)   \]  
where $\phi(a,T)$ is a certain integer depending only on $a$ and the genus of $T$.
\end{thm}

\section{Fourier coefficients of harmonic theta series}

In this section we rewrite the formula for the Fourier coefficients of a harmonic theta series as a sum over a double coset and give a recursive construction of the coset representatives.\\[-2mm]

Let $\theta_{L,h,k}$ be a harmonic Siegel theta series as described in the last section. Then the Fourier coefficient $a(T)$ of $\theta_{L,h,k}$ at a positive definite, even matrix $T \in M_n(\Z)$ is given by 
\[ 
a(T) = \sum_{x \in \Gamma_T} \det(Q(x,h))^k 
\]
where $\Gamma_T = \{ x \in L^n \, | \, Q(x)=T \}$. 

Right multiplication by elements in the group $\Or(T) = \{ \varepsilon \in \GL_n(\Z) \, | \, \varepsilon^T T \varepsilon = T \}$ defines an equivalence relation on $\Gamma_T$. The Fourier coefficient $a(T)$ vanishes unless $\det(\varepsilon)^k = 1$ for all $\varepsilon \in \Or(T)$ and 
\[  
a(T) = |\Or(T)| \sum_{x\Or(T) \in \Gamma_T/\Or(T)} \det(Q(x,h))^k 
\]
in this case. 

The group $\Or(L)$ acts on $\Gamma_T$ from the left by acting on each of the components. Let $\Or(L)_h \subset \Or(L)$ be the subgroup of elements leaving the span of $h_1, \ldots ,h_n$ in $V_{\C}$ invariant. Then for every $\sigma \in \Or(L)_h$ there exists a unique matrix $m_{\sigma} \in \GL_n(\C)$ such that
$( \sigma^{-1}(h_1),\ldots,\sigma^{-1}(h_n) ) 
= (h_1,\ldots,h_n) m_{\sigma}$
and 
\[  
Q(\sigma x,h) = Q(x,h) m_{\sigma}
\] 
for all $x \in \Gamma_T$. Note that $\theta_{L,h,k}$ vanishes if $\det(m_{\sigma})^k \neq 1$ for some $\sigma \in \Or(L)_h$. 

Assume that $\det(\varepsilon)^k = \det(m_{\sigma})^k = 1$ for all $\varepsilon \in \Or(T)$ and $\sigma \in \Or(L)_h$. Then $\det(Q(x,h))^k$ is constant on the double cosets $\Or(L)_h\backslash \Gamma_T/\Or(T)$ and we can write
\[
a(T) = \sum_{\Or(L)_h x \Or(T) \in \Or(L)_h\backslash \Gamma_T/\Or(T)} 
                                 |\Or(L)_h x \Or(T)| \det(Q(x,h))^k \, .  
\]

More generally, we have the following result.

\begin{lem}
Let $H$ be a subgroup of $\Or(L)$ containing $\Or(L)_h$ and $\mathcal{R}$ a set of representatives of $H \backslash \Gamma_T / \Or(T)$. Suppose $\det(m_{\sigma})^k = \det(\varepsilon)^k = 1$ for all $\sigma \in \Or(L)_h$ and $\varepsilon \in \Or(T)$. Then 
\[
a(T) = \frac{|\Or(L)_h|}{|H|} 
       \sum_{x \in \mathcal{R}}  |H x \Or(T)| 
       \sum_{\Or(L)_h \sigma \in \Or(L)_h \backslash H}\det(Q(\sigma x, h))^k \, .
\]
\end{lem}

\noindent
{\em Proof:} 
Since $\Gamma_T$ is invariant under $\Or(L)$, we can write 
\[
a(T) = \sum_{x \in \Gamma_T} \det(Q(x,h))^k = 
\frac{1}{|H|}\sum_{x \in \Gamma_T} \, \sum_{\sigma \in H}\det(Q(\sigma x, h))^k \, .
\]
The inner sum is constant on the double cosets $H \backslash \Gamma_T/ \Or(T)$ so that
\begin{align*}
a(T)
& = 
\frac{1}{|H|} \: \sum_{x\in \mathcal{R}} \; \sum_{\sigma \in H} 
                               \: |H x \Or(T)| \det(Q(\sigma x,h))^k \\
& =
\frac{1}{|H|} \: \sum_{x\in \mathcal{R}} \; \sum_{\sigma'\in \Or(L)_h} \; 
              \sum_{\Or(L)_h \sigma \in  \Or(L)_h\backslash H} 
              \: |H x \Or(T)| \det(Q(\sigma' \sigma x,h))^k \\
& = 
\frac{|\Or(L)_h|}{|H|} \: \sum_{x\in \mathcal{R}} \; 
                       \sum_{\Or(L)_h\sigma \in \Or(L)_h\backslash H}
              \: |H x \Or(T)| \det(Q(\sigma x,h))^k \, .
\end{align*}
Here we used $Q(\sigma'\sigma x,h) = Q(\sigma x,h) m_{\sigma '}$ and $\det(m_{\sigma '})^k = 1$. \eop\\[-1.5mm]

This formula allows the calculation of the Fourier coefficients of harmonic theta series in many interesting cases. 

A set of representatives $\mathcal{R}$ for $H \backslash \Gamma_T / \Or(T)$ can be constructed as follows. Let $T_m$ be the upper left $m \! \times \! m$-submatrix of $T$ and
\[
\Gamma_T^{m} = \{ x^m \in L^m \, | \, Q(x^m) = T_m  \} \, .
\]
Note that $\Gamma_T^n = \Gamma_T$. 

\begin{lem} \label{constS}
Let $2 \leq m \leq n$ and $\mathcal{S}_{m-1}$ be a set of representatives for $H\backslash \Gamma^{m-1}_T$. Let 
\[
\Gamma_{T,\mathcal{S}_{m-1}}^m = 
\{ x^m = (x^{m-1},x') \in \mathcal{S}_{m-1}\times L \, | \, Q(x^m) =  T_m \} 
                             \subset \Gamma_T^m
\]
and define an equivalence relation $\sim$ on $\Gamma_{T,\mathcal{S}_{m-1}}^m$ by $(x^{m-1},x')\sim (y^{m-1},y')$ if $x^{m-1}=y^{m-1}$ and there is some $\sigma\in H$ such that $\sigma x^{m-1}=x^{m-1}$ and $\sigma x'=y'$. Then a set of representatives for the equivalence classes of $\sim$ is a set of representatives for $H\backslash \Gamma_T^m$.
\end{lem}

\noindent
{\em Proof:}
Let $\mathcal{S}'_m$ be a set of representatives for $\Gamma_{T,\mathcal{S}_{m-1}}^m/\!\sim$. First we show that for every $x \in \Gamma_T^m$ there is an element $\sigma\in H$ such that $\sigma x\in \mathcal{S}'_m$. Write $x=(x^{m-1},x')$ with $x^{m-1}\in \Gamma_T^{m-1}$ and $x'\in L$. Since $\mathcal{S}_{m-1}$ is a set of representatives for $H \backslash \Gamma_T^{m-1}$, there is an element $\sigma\in H$ such that $\sigma x^{m-1}\in \mathcal{S}_{m-1}$. Then $\sigma x=(\sigma x^{m-1},\sigma x')$ is an element of $\Gamma_{T,\mathcal{S}_{m-1}}^m$ and therefore there is an element $\tau\in H$ such that $\tau\sigma x\in \mathcal{S}'_m$. Next, let $x,y\in \mathcal{S}'_m$ satisfy $\sigma x = y$ for some $\sigma\in H$. Write $x=(x^{m-1},x')$ and $y=(y^{m-1},y')$ with $x^{m-1},y^{m-1}\in \mathcal{S}_{m-1}$ and $x',y'\in L$. Since $\mathcal{S}_{m-1}$ is a set of representatives for $H\backslash \Gamma_T^{m-1}$ and $\sigma x^{m-1}=y^{m-1}$, we must have $x^{m-1}=y^{m-1}$ and therefore $\sigma x^{m-1}=x^{m-1}$. It follows that $x \sim y$, so $x=y$ as $\mathcal{S}'_{m}$ is a set of representatives for $\Gamma_{T,\mathcal{S}_{m-1}}^m/\!\sim$. \eop\\[-1.5mm]

Using the lemma we can construct a set $\mathcal{S}=\mathcal{S}_n$ of representatives for $H \backslash \Gamma_T$ starting from a set $\mathcal{S}_1$ of representatives of $H \backslash \Gamma_T^1$.

The set $\Gamma_{T,\mathcal{S}_{m-1}}^m$ is a subset of the finite set $\mathcal{S}_{m-1}\times L_{t_{m}}$ where $t_{m}$ is the lower right entry of $T_m$ and $L_{t_{m}}$ the set of norm $t_{m}$-vectors in $L$. So in order to construct $\mathcal{S}_{m}$ we can first determine the finite set $\mathcal{S}_{m-1}\times L_{t_{m}}$, then compute the subset $\Gamma_{T,\mathcal{S}_{m-1}}^m$ and finally a set of representatives for $\Gamma_{T,\mathcal{S}_{m-1}}^m/\!\sim$.

We put $\Gamma_{T,\mathcal{S}_0}^1 = \Gamma_T^1$ and let $x = (x_1,\ldots,x_n) \in \Gamma_T$ with $x_j \in L$. Define $x^1 = x_1 \in \Gamma_{T,\mathcal{S}_0}^1$. Then $\sigma_{x^1}x^1 \in S_1$ for some $\sigma_{x^1} \in H$. For $2 \leq j \leq n$ we recursively construct $x^j = \sigma_{x^{j-1}}\ldots\sigma_{x^1} (x_1,\ldots,x_j) \in \Gamma_{T,\mathcal{S}_{j-1}}^j$. The proof of Lemma \ref{constS} shows that there is some element $\sigma_{x^j} \in \Stab_H(\sigma_{x^{j-1}}x^{j-1})$ such that $\sigma_{x^j}x^j \in \mathcal{S}_{j}$. Then the element $\sigma_x = \sigma_{x^n}\ldots\sigma_{x^1} \in H$ satisfies $\sigma_x x \in \mathcal{S}$. It follows

\begin{lem} \label{constR}
Two elements $x,y \in \mathcal{S}$ are equivalent modulo $\Or(T)$ if and only if there is an element $\varepsilon \in \Or(T)$ such that 
%$\sigma_{x\varepsilon} x\varepsilon \in Hy$.
$\sigma_{x\varepsilon} x\varepsilon = y$.
\end{lem}

By means of Lemma \ref{constS} and \ref{constR} we can now easily construct a set $\mathcal{R}$ of double coset representatives for $H \backslash \Gamma_T / \Or(T)$. 

Let $x \in \mathcal{R}$ and $Hx\Or(T) = \bigcup_i H x \varepsilon_i$ be a disjoint decomposition of $Hx\Or(T)$ into cosets of $H$. Then
\[  |Hx\Or(T)| = \sum_i |Hx\varepsilon_i|  \, . \]

\section{A basis of $\S_{14}(\Sp_{12}(\Z))$} \label{basissec}

In this section we construct nine harmonic theta series of degree $6$ and weight $14$ associated with Niemeier lattices. We show that they are linearly independent by computing sufficiently many Fourier coefficients using the results from the previous section. It follows that these functions form a basis of $\S_{14}(\Sp_{12}(\Z))$.

We refer to \cite{CS} for the results on Niemeier lattices that we use. 
\\[-2mm]

The simplest example of a cusp form of degree $6$ and weight $14$ is the Ikeda lift 
%(cf.\ Theorem \ref{Ikedalift}) 
$F$ of the unique normalized cusp form $f = E_{10} \Delta \in \S_{22}(\SL_2(\Z))$. The function $f$ corresponds to the cusp form 
\begin{align*}
g(\tau) 
&= \frac{\eta(2\tau)^{19}\eta(4\tau)^{10}}{\eta(\tau)^6} + 4 \frac{\eta(\tau)^2\eta(4\tau)^{26}}{\eta(2\tau)^5}  \\[1mm]
&= q^3 + 10q^4 - 88q^7 - 132q^8 + 1275q^{11} + 736q^{12} - 8040q^{15} - 2880q^{16} + \ldots
\end{align*}
% q^3*eta(q^2)^19*eta(q^4)^10/eta(q)^6 + 4*q^4*eta(q)^2*eta(q^4)^26/eta(q^2)^5
% in $\S^+_{23/2}(\tilde{\Gamma}_0(4))$ 
under the Shimura correspondence. The Fourier coefficients of $F$ can be determined with the formula in Theorem \ref{Ikedalift}. We list a few of them below. 

%\vspace*{2mm}
\[
\renewcommand{\arraystretch}{1.2}
\begin{array}{c|c|c|c|c|c|c|c}
 T   & A_1^6     & A_2^3   & A_3^2 & A_6  & D_6 & E_6 & E_6'(3) \\[0.5mm] \hline
     &          &         &       &     &     &     & \\[-4.5mm]
a(T) & 63744000 & -128844 & 17600 & -88 & 10  & 1   & 100283601960
\end{array}
\]

\vspace*{1mm}
\noindent
More coefficients are given in Theorem \ref{coefficients}. They will provide a consistency check of our result on the Fourier coefficients of the harmonic theta series that we consider.

One of the basis elements will be a harmonic theta series for the Leech lattice. We describe its construction in detail.

Let $\mathcal{G}$ be the binary Golay code, i.e. the linear subspace of $\F_2^{24}$ spanned by the columns of 
\[
\left(
\begin{array}{cccccccccccc|cccccccccccc}
1 & 0 & 0 & 0 & 0 & 0 & 0 & 0 & 0 & 0 & 0 & 0 & 0 & 0 & 1 & 1 & 0 & 0 & 1 & 1 & 1 & 1 & 0 & 1 \\
0 & 1 & 0 & 0 & 0 & 0 & 0 & 0 & 0 & 0 & 0 & 0 & 1 & 0 & 0 & 1 & 1 & 0 & 0 & 1 & 1 & 1 & 1 & 0 \\
0 & 0 & 1 & 0 & 0 & 0 & 0 & 0 & 0 & 0 & 0 & 0 & 1 & 1 & 0 & 0 & 1 & 1 & 0 & 0 & 0 & 1 & 1 & 1 \\
0 & 0 & 0 & 1 & 0 & 0 & 0 & 0 & 0 & 0 & 0 & 0 & 0 & 1 & 1 & 0 & 0 & 1 & 1 & 0 & 1 & 0 & 1 & 1 \\
0 & 0 & 0 & 0 & 1 & 0 & 0 & 0 & 0 & 0 & 0 & 0 & 0 & 1 & 1 & 1 & 1 & 0 & 0 & 1 & 0 & 0 & 1 & 1 \\
0 & 0 & 0 & 0 & 0 & 1 & 0 & 0 & 0 & 0 & 0 & 0 & 1 & 0 & 1 & 1 & 1 & 1 & 0 & 0 & 1 & 0 & 0 & 1 \\
0 & 0 & 0 & 0 & 0 & 0 & 1 & 0 & 0 & 0 & 0 & 0 & 1 & 1 & 0 & 1 & 0 & 1 & 1 & 0 & 1 & 1 & 0 & 0 \\
0 & 0 & 0 & 0 & 0 & 0 & 0 & 1 & 0 & 0 & 0 & 0 & 1 & 1 & 1 & 0 & 0 & 0 & 1 & 1 & 0 & 1 & 1 & 0 \\
0 & 0 & 0 & 0 & 0 & 0 & 0 & 0 & 1 & 0 & 0 & 0 & 1 & 0 & 0 & 1 & 0 & 1 & 1 & 1 & 0 & 0 & 1 & 1 \\
0 & 0 & 0 & 0 & 0 & 0 & 0 & 0 & 0 & 1 & 0 & 0 & 1 & 1 & 0 & 0 & 1 & 0 & 1 & 1 & 1 & 0 & 0 & 1 \\
0 & 0 & 0 & 0 & 0 & 0 & 0 & 0 & 0 & 0 & 1 & 0 & 0 & 1 & 1 & 0 & 1 & 1 & 0 & 1 & 1 & 1 & 0 & 0 \\
0 & 0 & 0 & 0 & 0 & 0 & 0 & 0 & 0 & 0 & 0 & 1 & 0 & 0 & 1 & 1 & 1 & 1 & 1 & 0 & 0 & 1 & 1 & 0
\end{array}\right)^T.
\]
We define $V = \R^{24}$ and write the elements of $V$ as column vectors. We equip $V$ with the scalar product $(x,y) = \frac{1}{8} \sum_{i=1}^{24} x_i y_i$. Then the Leech lattice $\Lambda$ consists of the vectors
\begin{align*}
&0 + 2c + 4x , \\
&1 + 2c + 4y
\end{align*}
where $0=(0,\ldots,0)^T$, $1=(1,\ldots,1)^T$, $c \in \mathcal{G}$ (regarding the components of $c$ as real 0's and 1's) and $x,y \in \Z^{24}$ satisfy $\sum_{i=1}^{24} x_i = 0 \! \mod 2$ and  $\sum_{i=1}^{24} y_i = 1 \! \mod 2$.

We define $h=(h_1,\ldots,h_6) \in V_{\C}^6$ by
\[ h = 
\left(\begin{array}{cccccccccccccccccccccccc}
i & 0 & 0 & 0 & 0 & 0 & 0 & 0 & 1 & 0 & 0 & 0 & 0 & 0 & 0 & 0 & 0 & 0 & 0 & 0 & 0 & 0 & 0 & 0\\
0 & 0 & 0 & 1 & 0 & 0 & 0 & 0 & 0 & 0 & 0 & 0 & 0 & 0 & 0 & 0 & 0 & 0 & 0 & 0 & 0 & 0 & 0 & i\\
0 & 0 & 0 & 0 & 1 & 0 & 0 & 0 & 0 & 0 & 0 & 0 & 0 & 0 & 0 & 0 & 0 & 0 & 0 & 0 & i & 0 & 0 & 0\\
0 & 0 & 0 & 0 & 0 & 0 & 0 & 1 & 0 & 0 & 0 & 0 & 0 & 0 & 0 & 0 & 0 & 0 & i & 0 & 0 & 0 & 0 & 0\\
0 & 0 & 0 & 0 & 0 & 0 & 0 & 0 & 0 & 1 & 0 & 0 & 0 & 0 & 0 & i & 0 & 0 & 0 & 0 & 0 & 0 & 0 & 0\\
0 & 0 & 0 & 0 & 0 & 0 & 0 & 0 & 0 & 0 & 0 & 0 & 1 & 0 & 0 & 0 & 0 & 0 & 0 & i & 0 & 0 & 0 & 0
\end{array}\right)^T.
\]
Then 
\[ \theta_{\Lambda,h,2}(Z) = \sum_{x \in \Lambda^6} \det(Q(x,h))^2 e^{\pi i \tr(Q(x)Z) }  \]
is in $\S_{14}(\Sp_{12}(\Z))$.

We explain the choice of $h$. The automorphism group of the Golay code is the Mathieu group $M_{24}$. The stabilizer $\Aut(\mathcal{G})_c$ of the dodecad 
\[
c = (1, 0, 0, 1, 1, 0, 0, 1, 1, 1, 0, 0, 1, 0, 0, 1, 0, 0, 1, 1, 1, 0, 0, 1)^T \in \mathcal{G}
\]
in $\Aut(\mathcal{G})$ is isomorphic to $M_{12}$. Let $e_j=(0,\ldots,0,1,0,\ldots,0)^T \in V$ with $1$ at the $j$-th position and $U$ the span of the $e_j$, $j\in \text{supp}(c)$. Then the only element in $\Aut(\mathcal{G})_c$ acting trivially on $U$ is the identity because $M_{12}$ is simple. The permutation 
\[
\tau = (1, 9)(2, 12)(3, 7)(4, 24)(5, 21)(6, 18)(8,19)(10, 16)(11, 17)(13, 20)(14, 23)(15, 22)
\]
has cycle type $2^{12}$ and is in $\Aut(\mathcal{G})_c$. Its centralizer in $\Aut(\mathcal{G})_c$ is a maximal subgroup of order $240$ which we denote by $C_{c,\tau}$. The vectors $h_i$ are supported on those cycles $(j,k)$ of $\tau$ for which $c_j$ and hence also $c_k$ is equal to $1$, i.e. on the cycles 
\[
(1, 9),(4, 24),(5, 21),(8,19),(10, 16),(13, 20) \, . 
\]
Since the elements of $\Aut(\mathcal{G})$ are permutations of coordinates, this group embeds naturally into $\Or(\Lambda)$. Furthermore, $\mathcal{G}$ acts on $\Lambda$ by sign changes on the support of an element. Recall that $\Or(\Lambda)_h$ is the subgroup of $\Or(\Lambda)$ which preserves the complex span of the $h_i$.

\begin{lem}
The group $\Or(\Lambda)_h$ has order $30720$ and is contained in $\mathcal{G}C_{c,\tau} \subset \Or(\Lambda)$. Moreover, $\det(m_{\sigma})= \pm 1$ for all $\sigma \in \Or(\Lambda)_h$.
\end{lem}

\noindent
{\em Proof:}
Let $\vartheta$ be the element in $\Or(\Lambda)_h$ corresponding to $c$. Then the eigenspaces of $\vartheta$ are $U$ and $U^{\perp}$. An element $\sigma \in \Or(\Lambda)_h$ preserves the space spanned by the vectors $h_1,\ldots,h_6$ and the complex conjugates $\ov{h}_1,\ldots,\ov{h}_6$ and hence also the spaces $U$ and $U^{\perp}$. Therefore $\sigma$ commutes with $\vartheta$. Since the centralizer of $\vartheta$ in $\Or(\Lambda)$ is $\mathcal{G}\Aut(\mathcal{G})_c$ (see \cite{CS}, chapter 10, section 3.7), we can write $\sigma = \sigma_1\sigma_2$ with $\sigma_1 \in \mathcal{G}$ and $\sigma_2 \in \Aut(\mathcal{G})_c$. We first consider $h_1 = ie_1 + e_9$. Then $\sigma(h_1) = \sigma_1(ie_{\sigma_2(1)} + e_{\sigma_2(9)})=\pm ie_{\sigma_2(1)} \pm e_{\sigma_2(9)}$ where the signs need not be the same. Since $\sigma(h_1)$ is in the space generated by the vectors $h_1,\ldots,h_6$, we have $\tau(\sigma_2(1)) = \sigma_2(9) = \sigma_2(\tau(1))$. A similar argument for the other $h_j$ shows that the restriction of $\tau\sigma_2$ to $U$ is equal to the restriction of $\sigma_2\tau$ to $U$. It follows that $\tau\sigma_2 = \sigma_2\tau$. This implies that $\Or(\Lambda)_h$ is contained in $\mathcal{G}C_{c,\tau}$. Now we can determine $\Or(\Lambda)_h$ explicitly with a computer and verify the statements in the lemma. \eop\\[-1.5mm]

We use the results from the previous section to determine some Fourier coefficients of $\theta_{\Lambda,h,2}$. In our computations we have chosen $H$ to be the centralizer of the element $\vartheta$ corresponding to $c$ in $\Or(\Lambda)$, i.e.  $H = \mathcal{G}\Aut(\mathcal{G})_{c}$, so the index of $\Or(\Lambda)_h$ in $H$ is $12672$.

Similarly, we calculate Fourier coefficients of harmonic theta series associated with the Niemeier lattices $N(A_6^4), N(D_6^4), N(E_6^4), N(A_4^6), N(D_4^6), N(A_3^8), N(A_2^{12})$ and $N(A_1^{24})$. 
The realizations of these lattices and the choices of $h$ can be found in the appendix. We obtain the following theorem.

\begin{thm} \label{coefficients}
The Fourier coefficients of the normalized harmonic theta series and the Ikeda lift $F$ of $E_{10}\Delta$ are given by
\[
\begin{array}{c|rrrrr}
    & N(A_6^4) & N(D_6^4) & N(E_6^4) & N(A_4^6) & N(D_4^6) \\[0.5mm] \hline 
    & & & & & \\[-3.5mm]
A_6 & 1 & 0 & 0 & 0 & 0  \\
D_6 & 0 & 1 & 0 & 0 & 0  \\
E_6 & 0 & 0 & 1 & 0 & 0  \\
A_5A_1 & -12 & -30  & -20  & 0 & 0  \\
D_5A_1 & 0 & -10  & -32 & 0 & 0  \\
A_4A_2 & 30 & 120  & 240 & 1  & 0  \\
D_4A_2 & 0 & 72  & 192 & 0 & 1  \\
A_4A_1^2 & 40 & 260  & 320 & -2  & 0 \\
D_4A_1^2 & 0 & 48  & 384 & 0 & -2  \\
A_3^2 & -32 & -192  & -432 & 0  & 0  \\
A_3A_2A_1 & -96 & -276  & -480 & -8 & -6 \\
A_3A_1^3 & 576  & 672  & 192 & 48  & 36  \\
A_2^3 & 648 & 1296  & 900 & 54 & 36  \\
A_2^2A_1^2 & -432  & -1080  & -1152 & -36 & 0  \\
A_2A_1^4 & -1152 & -3456  & 768 & 96 & -152  \\
A_1^6 & -11520 & 3840  & -46080 & -2880 & 240   \\
A_1(2)A_5 & -816 & -11436 & -11568 & 0 & 0  \\
A_1(2)D_5 & 0 & -1376  & -9600 & 0 & 0  \\
A_1(2)A_4A_1 & -2320 & 61020  & 121440 & -4  & 0  \\
A_1(2)D_4A_1 & 0 & -2784  & 88320 & 0 & -12   \\
A_1(2)A_3A_2 & -27072 & -94464  & -156384 & -288 & -792  \\
A_1(2)A_3A_1^2 & 127744 & 436736  &627456 & 320 & 1216  \\
A_1(2)A_2^2A_1 & -2592 & -282528  & -1221696 & -4104 & 1512 \\
A_1(2)A_2A_1^3 & -208512 & -1416384 & -743040  & 24480  & 1656   \\
A_1(2)A_1^5 & -2772480 & -5598720  & -19937280 & -191360  & -84320   \\
E_6(2) & 18247680 & -114048000 & 436423680 & -2142720 & 11986560 \\
E_6'(3) & 149532480 & -874800000 & 2327826600 & -15843600 & 44621280 \\[0.5mm] \hline 
    & & & & & \\[-3.5mm]
        & -88 & 10 & 1 & -6840 & 1872   
\end{array}
\]
\[
\begin{array}{c|rrrrr}
    &  N(A_3^8) & N(A_2^{12}) & N(A_1^{24}) & \Lambda & F\\[0.5mm] \hline 
    & & & & & \\[-3.5mm]
A_6  & 0 & 0 & 0 & 0 & -88 \\
D_6  & 0 & 0 & 0 & 0 & 10\\
E_6  & 0 & 0 & 0 & 0 & 1\\
A_5A_1  & 0 & 0 & 0 & 0 & 736\\
D_5A_1  & 0 & 0 & 0 & 0 & -132\\
A_4A_2  & 0 & 0 & 0 & 0 & -8040 \\
D_4A_2  & 0 & 0 & 0 & 0 & 2784\\
A_4A_1^2  & 0 & 0 & 0 & 0 & 13080\\
D_4A_1^2 & 0 & 0 & 0 & 0 & -2880\\
A_3^2  & 1  & 0 & 0 & 0 & 17600\\
A_3A_2A_1 & -6  & 0 & 0 & 0 & -54120\\
A_3A_1^3  & 20 & 0 & 0 & 0 & 38016\\
A_2^3  & 0  & 1 & 0 & 0 & -128844\\
A_2^2A_1^2  & 72 & -2   & 0 & 0 & 1073520 \\
A_2A_1^4  & -480 & 12 & 0 & 0 & -6503424\\
A_1^6   & 4000  & -120  & 0 & 0 & 63744000\\
A_1(2)A_5  & 0 & 0 & 0 & 0 & -54120\\
A_1(2)D_5  & 0 & 0 & 0 & 0 & -23360 \\
A_1(2)A_4A_1   & 0 & 0 & 0 & 0 & 963160 \\
A_1(2)D_4A_1  & 0 & 0 & 0 & 0  & 38016 \\
A_1(2)A_3A_2  & -150 & 0 & 0 & 0 & -801792\\
A_1(2)A_3A_1^2 & -816 & 0 & 0 & 0  & -20142080\\
A_1(2)A_2^2A_1  & 1080 & 12  & 0 & 0 & 48185280\\
A_1(2)A_2A_1^3  & 11208  & -72 & 0 & 0 & 15586560 \\
A_1(2)A_1^5  & -98880  & 560  & 4 & 0  & -841420800\\
E_6(2)  & -18809280 & 5212800 & -27081 & 1 & 47850946560\\
E_6'(3)  & -90357120 & 21107880 & -108864 & 4 & 100283601960\\[0.5mm] \hline 
    & & & & & \\[-3.5mm]
     & 17136 & 216288 & -146810880 & -4767869952000
\end{array}
\]
The last column is a linear combination of the preceding ones. The corresponding coefficients are given in the last row.
\end{thm}

\noindent
Ta\"ibi showed that the dimension of $S_{14}(\Sp_{12}(\Z))$ is $9$ (\cite{dimformula}, Section 5.5., Table 3). This implies

\begin{thm} \label{basis}
Any nine of the above functions form a basis of $S_{14}(\Sp_{12}(\Z))$. In particular the nine harmonic theta series span $S_{14}(\Sp_{12}(\Z))$.
\end{thm}

\noindent
Recall that Proposition \ref{both} does not give the existence of a basis consisting of harmonic theta series because the condition $m/2>2n$ is not satisfied for $m=24$ and $n=6$.

\section{A lower bound for the Kodaira dimension of $\mathcal{A}_6$}

In this section we show that the Kodaira dimension of $\A_6$ is non-negative.\\[-2mm]

Let $X$ be a complex smooth projective variety of dimension $n$. For a non-negative integer $d$ the $d$-th plurigenus $P_d$ of $X$ is defined as the dimension of the space of global sections of $\omega^d$ where $\omega$ is the canonical bundle of $X$, i.e. $P_d = \dim (H^0(X,\omega^d))$. The Kodaira dimension $k(X)$ of $X$ is defined as $-\infty$ if $P_d=0$ for all $d>0$. Otherwise, $k(X)$ is the smallest non-negative integer $k$ such that $P_d(X)/d^k$ is bounded. We have $k(X) \leq \dim(X)$ and say that $X$ is of general type if $k(X) = \dim(X)$. Furthermore, $k(X)$ is a birational invariant. The Kodaira dimension of a variety $Y$ is defined as the Kodaira dimension of a smooth projective variety birational to $Y$.

The quotient $\A_n = \Sp_{2n}(\Z) \backslash H_n$ parametrizes the principally polarized complex abelian varieties of dimension $n$. It has the structure of a normal quasi-projective variety of dimension $n(n+1)/2$ and is called the Siegel modular variety of degree $n$. The Kodaira dimension of $\A_n$ for $n \neq 6$ has been known for more than 30 years. Namely $\A_n$ is of general type for $n \geq 7$ and $k(\A_n) = - \infty$ for $n \leq 5$. 

We can use Siegel modular forms to construct global sections of $\omega^d$. Let $n \geq 3$. We write $Z \in H_n$ as
\[ 
Z = 
\begin{pmatrix}
z_{11}  & \cdots & z_{1n}  \\
\vdots & \ddots & \vdots \\
z_{1n}  & \cdots & z_{nn}  
\end{pmatrix} \, . 
\]
Let $H^{\circ}_n \subset H_n$ be the subset of points whose stabilizer in $\Sp_{2n}(\Z)$ is $\{ \pm I_{2n} \}$. Then $H^{\circ}_n$ is an open subset of $H_n$ and $\A_n^{\circ} =  \Sp_{2n}(\Z) \backslash H^{\circ}_n$ is the smooth locus of $\A_n$ (see \cite{F}, Hilfssatz III.5.15). Let $f$ be a Siegel modular form of degree $n$ and weight $(n+1)k$. Then the transformation property of $f$ implies that
\[  
f(Z)(dz_{11} \wedge dz_{12} \wedge \ldots \wedge dz_{nn})^k
\]
is a global section of $\omega_{\A_n^{\circ}}^k$. Tai has shown that (cf. \cite{Tai} and \cite{F}, Satz III.5.24)

\begin{prp} \label{taissatz}
Suppose $n \geq 5$ and the vanishing order of $f$ at $\infty$ is at least $k$. Let $\ov{\A_n}$ be a smooth compactification of $\A_n$. Then 
\[  
f(Z)(dz_{11} \wedge dz_{12} \wedge \ldots \wedge dz_{nn})^k
\]
can be extended to a global section of $\omega_{\ov{\A_n}}^k$.
\end{prp}

Let $\theta_{\Lambda,h,2}$ be the harmonic theta series from the previous section. It is non-zero, has weight $14$ and vanishing order $2$ at $\infty$ because the Leech lattice has no elements of norm $2$. Hence it defines a non-trivial section of $\omega_{\ov{\A_6}}^2$ by Proposition \ref{taissatz}. We obtain

\begin{thm} \label{kodairadim}
The Kodaira dimension of $\A_6$ is non-negative.
\end{thm}

\noindent
According to Table 6 in \cite{CT} the space $\S_{7}(\Sp_{12}(\Z))$ is trivial which means that the canonical bundle $\omega_{\ov{\A_6}}$ has no non-trivial sections. We showed that the bicanonical bundle $\omega_{\ov{\A_6}}^2$ admits a non-trivial section. 

A consequence of the above theorem is

\begin{cor}
$\A_n$ is unirational if and only if $n \leq 5$.
\end{cor}

\section{Appendix}

In Section \ref{basissec} we constructed a basis of $\S_{14}(\Sp_{12}(\Z))$ consisting of harmonic theta series $\theta_{N,h,2}$. Here we list the realizations of the Niemeier lattices $N$ different from the Leech lattice and the $h=(h_1,\ldots,h_6)$ that we used.\\[-2mm] 

Let $V = \R^{24}$. As before we write the elements of $V$ as column vectors. 
%\\[-2mm] 

We begin with the Niemeier lattices of type $N(A_n^j)$. We denote by $G(A_n)$ the standard Gram matrix of $A_n$.

To define $N(A_1^{24})$ we equip $V$ with inner product given by twice the standard inner product. For every vector $c = (c_1, \ldots, c_{24})^T$ in the Golay code $\mathcal{G}$ define the vector $y^c = \frac{1}{2} c \in \R^{24}$ (regarding the components of $c$ as real 0's and 1's). Then $N(A_1^{24})$ is the lattice in $V$ generated by $\Z^{24}$ and the vectors $(y^c)_{c \in \mathcal{G}}$.

For $N(A_6^4)$ we equip $V$ with inner product given by the Gram matrix $G(A_6)^4$, i.e. with the block diagonal matrix with four blocks, each equal to $G(A_6)$. We define the vector $v = \frac{1}{7}(1,2,3,4,5,6)^T \in \R^6$. Then $N(A_6^4)$ is the lattice in $V$ generated by $\Z^{24}$ and the vectors 
\[  (v,2v,v,6v)^T, (v,v,6v,2v)^T \, . \]

To define $N(A_4^6)$ we equip $V$ with inner product given by the Gram matrix $G(A_4)^6$. Let $v = \frac{1}{5}(2 , 4, -4, -2)^T  \in \R^4$. Then $N(A_4^6)$ is the lattice generated by $\Z^{24}$ and the vectors 
\[  (v , 0 , v , 4v , 4v , v)^T, (v ,v,4v ,4v ,v ,0)^T, (v ,4v , 4v ,v ,0  ,v)^T \, . \]

For $N(A_3^8)$ we equip $V$ with inner product given by the Gram matrix $G(A_3)^8$. Define $v = \frac{1}{4}(1 , 2, -1)^T \in \R^3$. Then $N(A_3^8)$ is the lattice generated by $\Z^{24}$ and the vectors 
\begin{align*}
(v, 3v , v  , 2v , v , 0  , 0  , 0)^T, 
(v, v  , 2v , v  , 0 , 0  , 0  , 3v)^T,
(v, 2v , v  , 0  , 0 , 0  , 3v , v )^T,  
(v, v  , 0  , 0  , 0 , 3v , v  , 2v)^T.
\end{align*}

To define $N(A_2^{12})$ we equip $V$ with inner product given by the Gram matrix $G(A_2)^{12}$. Let $v = \frac{1}{3}(1,2)^T \in \R^2$. Then $N(A_2^{12})$ is the lattice generated by $\Z^{24}$ and 
\begin{align*}
& 
(v , 0 , 0 , 0 , 0 , 0 , 0 , v , v , v ,v ,v )^T, 
(0 , v , 0 , 0 , 0 , 0 , v , 0 , v , 2v,2v ,v )^T, 
(0 , 0 , v , 0 , 0 , 0 , v , v , 0 , v ,2v ,2v )^T, \\
& 
(0 , 0 , 0 , v , 0 , 0 , v , 2v, v , 0 ,v ,2v )^T,
(0 , 0 , 0 , 0 , v , 0 , v , 2v, 2v, v ,0 ,v )^T,
(0 , 0 , 0 , 0 , 0 , v , v , v , 2v, 2v,v ,0 )^T.
\end{align*}

Next we describe our realization of the Niemeier lattices of type $N(D_n^j)$.

We equip $V$ with the standard inner product. We define the vectors $s=\frac{1}{2}(1 , 1 , 1 , 1 , 1 , 1)^T$, $v= (0 , 0 , 0 , 0 , 0 , 1)^T$ and $c=s-v$ in $\R^6$. Then $N(D_6^4)$ is the lattice generated by $\{ (x_1,\cdots,x_n)^T \in \Z^{24} \, | \, 
\sum_{i=1+6j}^{6+6j} x_i = 0 \mod 2 \, \, \text{ for } j = 0,\ldots,3 \}$ 
and the vectors 
\[  (0 , s ,v , c)^T, (0 , v ,c , s)^T, (s , 0 ,c , v)^T, (v , 0  , s , c)^T  \, . \]

The lattice $N(D_4^6)$ is realized analogously. We equip $V$ with the standard inner product and define $s=\frac{1}{2}(1 , 1 , 1 , 1)^T$, $v= (0 , 0 , 0 , 1)^T$ and $c=s-v$. Then $N(D_4^6)$ is the lattice generated by $\{ (x_1,\cdots,x_n)^T \in \Z^{24} \, | \, 
\sum_{i=1+4j}^{4+4j} x_i = 0 \mod 2 \, \, \text{ for } j = 0,\ldots,5 \}$ and the vectors
\[
(v , c , c , v , c , v)^T, (c ,v,v ,c ,c ,v)^T, (c ,v,c ,v ,v ,c)^T, (s ,0,c ,v ,0 ,s)^T, (0 ,0,c ,c ,c ,c)^T, (v ,0 , c ,s ,v  ,0)^T \, . 
\]

Finally let 
\[
G(E_6) =
\begin{pmatrix}
2 & -1 & 0 & 0& 0 & 0\\
-1& 2 & -1 & 0 & 0 & 0\\
0 & -1 & 2 &-1 & 0 & -1\\
0 & 0 & -1 & 2 &-1 & 0\\
0 & 0 & 0 & -1 &2 & 0\\
0 & 0 & -1 & 0 & 0 & 2
\end{pmatrix}
\]
be the Gram matrix of $E_6$. We equip $V$ with inner product given by the Gram matrix $G(E_6)^{4}$ and define the vector $v = \frac{1}{3} (1,-1,0,1,-1,0)^T\in \R^6$. Then $N(E_6^4)$ is the lattice in $V$ generated by $\Z^{24}$ and the vectors 
\[  (v,0,v,2v)^T, (v,2v,0,v)^T \, . \]

We define $h$ as follows.

For $N(A_1^{24})$ we let $h$ be as in the case of the Leech lattice.

For $N(A_6^4)$, $N(D_6^4)$ and $N(E_6^4)$ we put
\begin{align*}
h 
&= (h_1,\cdots,h_6)  \\
&= \begin{pmatrix}
1&0&0&0&0&0&i&0&0&0&0&0&0&0&0&0&0&0&0&0&0&0&0&0\\
0&1&0&0&0&0&0&i&0&0&0&0&0&0&0&0&0&0&0&0&0&0&0&0\\
0&0&1&0&0&0&0&0&i&0&0&0&0&0&0&0&0&0&0&0&0&0&0&0\\
0&0&0&1&0&0&0&0&0&i&0&0&0&0&0&0&0&0&0&0&0&0&0&0\\
0&0&0&0&1&0&0&0&0&0&i&0&0&0&0&0&0&0&0&0&0&0&0&0\\
0&0&0&0&0&1&0&0&0&0&0&i&0&0&0&0&0&0&0&0&0&0&0&0
\end{pmatrix}^T  \, . 
\end{align*}

For $N(A_4^6)$ and $N(D_4^6)$ we define 
\begin{align*}
h 
&= (h_1,\ldots,h_6)  \\
&= \begin{pmatrix}
1&0&0&0&0&0&0&0&i&0&0&0&0&0&0&0&0&0&0&0&0&0&0&0\\
0&1&0&0&0&0&0&0&0&i&0&0&0&0&0&0&0&0&0&0&0&0&0&0\\
0&0&1&0&0&0&0&0&0&0&i&0&0&0&0&0&0&0&0&0&0&0&0&0\\
0&0&0&1&0&0&0&0&0&0&0&i&0&0&0&0&0&0&0&0&0&0&0&0\\
0&0&0&0&1&0&0&0&0&0&0&0&i&0&0&0&0&0&0&0&0&0&0&0\\
0&0&0&0&0&1&0&0&0&0&0&0&0&i&0&0&0&0&0&0&0&0&0&0
\end{pmatrix}^T \, . 
\end{align*}

For $N(A_3^8)$ we choose 
\begin{align*}
h 
&= (h_1,\ldots,h_6)  \\
&= \begin{pmatrix}
1&0&0&0&0&0&i&0&0&0&0&0&0&0&0&0&0&0&0&0&0&0&0&0\\
0&1&0&0&0&0&0&i&0&0&0&0&0&0&0&0&0&0&0&0&0&0&0&0\\
1&2&3&0&0&0&0&0&0&0&0&0&0&0&0&0&0&0&0&0&0&0&0&i\sqrt{6}\\
0&0&0&1&0&0&0&0&0&i&0&0&0&0&0&0&0&0&0&0&0&0&0&0\\
0&0&0&0&1&0&0&0&0&0&i&0&0&0&0&0&0&0&0&0&0&0&0&0\\
0&0&0&1&2&3&0&0&0&0&0&0&0&0&0&0&0&0&0&0&i\sqrt{6}&0&0&0
\end{pmatrix}^T \, . 
\end{align*}

For $N(A_2^{12})$ we take 
\begin{align*}
h 
&= (h_1,\ldots,h_6)  \\
&= \begin{pmatrix}
1&0&0&0&0&0&0&0&0&0&0&0&0&0&0&0&i&0&0&0&0&0&0&0\\
0&1&0&0&0&0&0&0&0&0&0&0&0&0&0&0&0&i&0&0&0&0&0&0\\
0&0&1&0&0&0&0&0&0&0&0&0&i&0&0&0&0&0&0&0&0&0&0&0\\
0&0&0&1&0&0&0&0&0&0&0&0&0&i&0&0&0&0&0&0&0&0&0&0\\
0&0&0&0&1&0&0&0&0&0&0&0&0&0&0&0&0&0&i&0&0&0&0&0\\
0&0&0&0&0&1&0&0&0&0&0&0&0&0&0&0&0&0&0&i&0&0&0&0
\end{pmatrix}^T \, . 
\end{align*}


\begin{thebibliography}{mmmmx}

\addcontentsline{toc}{section}{References}

\bibitem[B]{Basisproblem} S.\ B\"ocherer, {\em \"Uber die Fourier-Jacobi-Entwicklung Siegelscher Eisensteinreihen. II}, Math.\ Z.\ {\bf 189} (1985), 81--110
\bibitem[C]{C} C.\ H.\ Clemens, {\em Double solids}, Adv.\ in Math.\ {\bf 47} (1983), 107--230 
\bibitem[CS]{CS} J.\ H.\ Conway, N.\ J.\ A.\ Sloane, {\em Sphere packings, lattices and groups}, Grundlehren der Math.\ Wiss.\ {\bf 290}, 3rd ed., Springer, New York, 1999
\bibitem[CT]{CT} G.\ Chenevier, O.\ Ta\"{i}bi, {\em Discrete series multiplicities for classical groups over $\Z$ and level 1 algebraic cusp forms}, Publ.\ Math.\ Inst.\ Hautes \'{E}tudes Sci.\ {\bf 131} (2020), 261--323
%arXiv:1907.08783 
\bibitem[D]{D} R.\ Donagi, {\em The unirationality of $\A_5$}, Ann.\ of Math.\ (2) {\bf 119} (1984), 269--307
\bibitem[F1]{F1} E.\ Freitag, {\em Der K\"orper der Siegelschen Modulfunktionen}, Abh.\ Math.\ Sem.\ Univ.\ Hamburg {\bf 47} (1978), 25--41
\bibitem[F2]{F2} E.\ Freitag, {\em Die Kodairadimension von K\"orpern automorpher Funktionen}, J.\ Reine Angew.\ Math.\ {\bf 296} (1977), 162--170 
\bibitem[F3]{F}  E.\ Freitag, {\em Siegelsche Modulfunktionen}, Grundlehren der Math.\ Wiss.\ {\bf 254}, Springer, Berlin, 1983
\bibitem[HM]{HM} J.\ Harris, D.\ Mumford, {\em On the Kodaira dimension of the moduli space of curves}, with an appendix by William Fulton, Invent.\ Math.\ {\bf 67} (1982), 23--88
\bibitem[I]{Ikedalift} T.\ Ikeda, {\em On the lifting of elliptic cusp forms to Siegel cusp forms of degree 2n}, Ann.\ of Math.\ (2) {\bf 154} (2001), 641--681 
\bibitem[K]{Likeda} W.\ Kohnen, {\em Lifting modular forms of half-integral weight to Siegel modular forms of even genus}, Math.\ Ann.\ {\bf 322} (2002), 787--809 
\bibitem[MM]{MoMu} S.\ Mori, S.\ Mukai, {\em The uniruledness of the moduli space of curves of genus 11}, Algebraic geometry (Tokyo/Kyoto, 1982), 334--353, Lecture Notes in Math.\ {\bf 1016}, Springer, Berlin, 1983
\bibitem[M]{Mum} D.\ Mumford, {\em On the Kodaira dimension of the Siegel modular variety}, Algebraic geometry--Open problems (Ravello, 1982), 348--375, Lecture Notes in Math.\ {\bf 997}, Springer, Berlin, 1983
\bibitem[Tai]{Tai} Y.-S. Tai, {\em On the Kodaira dimension of the moduli space of abelian varieties}, Invent.\ Math.\ {\bf 68} (1982), 425--439 
\bibitem[T]{dimformula} O.\ Ta\"{i}bi, {\em Dimensions of spaces of level one automorphic forms for split classical groups using the trace formula}, Ann.\ Sci.\ \'{E}c. Norm.\ Sup\'{e}r.\ 50 (2017), 269--344
\bibitem[V]{V} A.\ Verra, {\em A short proof of the unirationality of $\A_5$}, Nederl.\ Akad.\ Wetensch.\ Indag.\ Math.\ {\bf 46} (1984), 339--355


\end{thebibliography}
\end{document}